\documentclass[a4paper,11pt]{article}
\usepackage{amsmath, amsthm, amssymb, amsfonts}
\usepackage{mathtools}
\usepackage{graphicx}
\usepackage{caption}
\usepackage{subcaption}
\usepackage{hyperref}
\usepackage[utf8]{inputenc}
\usepackage[T1]{fontenc}
\usepackage{hyperref}
\usepackage{scrextend}
\usepackage{authblk}
\usepackage{float}
\usepackage{tikz-cd}
\usepackage[font=small,skip=0pt]{caption}
\deffootnote{2em}{0em}{\thefootnotemark\quad}

\newtheorem{prop}{Proposition}
\newtheorem{corollary}{Corollary}

\begin{document}
\title{ On the intimate association between even binary palindromic words and the Collatz-Hailstone iterations}
\author{T. Raptis \footnote{I.P.I., Portsmouth, UK}  \footnote{traptis@protonmail.com}} 
\date{}
\maketitle

\abstract{The celebrated $3x+1$ problem is reformulated via the use of an analytic expression of the trailing zeros sequence resulting in a single
branch formula $f(x)+1$ with a unique fixed point. The resultant formula $f(x)$ is also found to coincide with that of the discrete
derivative of the sorted sequence of fixed points of the reflection operator on even binary palindromes of fixed even length \textit{2k} 
in any interval $[0\cdots2^{2k}-1]$. A set of equivalent reformulations of the problem are also presented.}

\section{Introduction}

The so called, Hailstone sequences are the iterates of a discrete branched dynamical system over positive integers that were first introduced by Lothar Collatz \cite{col} at 1937. It became alternatively known as the \textit{'3x+1'} problem due to the particular formula used in one of its branches.

Lately, the notion of a generic \textit{'3x+1'} semigroup was also introduced \cite{lag}\cite{far}. At 1972, Conway had already
generalized the problem in the generic form $(a_i x+b_i), 0(mod b_i)$ in a manner that allowed him to create a programming language called
\textit{'fractran'} which was of the same power with a universal Turing machine thus being able to derive standard undecidability results\cite{con}.

A very recent review of the problem has been given by Lagarias\cite{rev}. Also, recent results on connections of this problem with cellular automata
appeared which are reviewed in \cite{ca} together with previous similar attemtps. Mostly, strong associations with Wang tiling machines have been introduced in the work of Sterin\cite{wang}
with possible connections to biological complexity.

An analysis of the resulting Hailstone sequences of iterants from a
physicist's perspective appeared in \cite{luz1} and \cite{luz2}.
It is hoped that the preset analysis will also be of interest for 
other physically inspired toy models for stochastic and fractal
processes.

In the next section, an appropriate reformulation of the Hailstone
iteration is introduced which makes use of a special function also
associated with the so called \textit{dyadic valuation}\cite{val}.
This allows transforming the original branched process to a single
branch one with a unique fixed point.

In section 3, a set of necessary definitions are introduced for
palindromic words or palindromes based on fixed maximal length
binary expansions as fixed points of the reflection group over
such expansions. A hierarchical construct is used to extract 
certain scaling maps associating each expansion length $L$ with
its next one across different intervals of exponential length
revealing the underlying tree structure of such patterns.

Certain properties of these hierarchies of patterns are discussed
and a crucial property is proven that provides a direct link 
with the original Hailstone process.

In section 4, the role of palindromes and the associated reflection
group is discussed revealing an interesting type of interaction
between a form of 'mirror' images inside the main process.

Furthermore, two indices in the form of binary probability measures
are proposed for the study of the conjectured global convergence,
associated with both the inner reflective structure as well as
the internal complexity of the binary patterns produced by the 
Hailstone process.

\section{ Reformulation of the Collatz-Hailstone (CH) iteration }

The standard, or $3x+1$ Collatz-Hailstone process is defined via the branched map

\begin{equation}
x_{n+1} = \left\{
\begin{array}{lr}
x_n/2, & 0(mod 2) \\
3x_n+1, & 1(mod 2)
\end{array}
\right.
\end{equation}

Let us introduce the two auxilliary maps $f_0(x) = x/2$ and $f_1(x) = 3x+1$. It is obvious that any invocation of $f_0$ would cause this iteration to enter a cycle any time $x_n$ reaches a power of two since it would then remain on the first branch until it reaches 1 with the second branch 
mapping $1 \rightarrow 4$ immediately after. 

On the other hand, $f_1$ always maps odd integers to even integers thus any iteration stays at the second branch only once while looping
over the first branch until all powers of two divisors are exausted.

The particular sequence of all integers after removal of its 2 factors
is already known as the odd part sequence, catalogued as \textit{A000265} in the OEIS database\cite{oeis265}

The resulting symbolic dynamics of branch execution then is intimately related with binary divisibility or equivalently, 
the amount of zeros present at the start of any binary expansion. This is another well known sequence
in computer science under the name of \textit{trailing zeros}(TZS), also catalogued as \textit{A007814} \cite{oeis}. 
Any integer is then represented as $x = \sigma_{odd}(x)2^{t(x)}$
where $t(x)$ denotes specific values of the TZS.

In terms of the run length analysis of symbolic sequences\cite{rl},\cite{rla} where every bit string of length $L$ is represnted by an alternating polynomial, the TZS corresponds to the zero order coefficient for all compressed binary expansions. This will be presented in more detail in the next sections.

When represented as a sequence over all integers, the TZS is equivalent  to the so called, \textit{2-adic valuation}, the first of the \textit{$\pi$-adic} 
valuations corresponding to the expression of all the exponents of prime factorizations as sequences\cite{piadic}. The structure of TZS encodes a special tree graph which in the context of word combinatorics is related to the \textit{Zimin words} or, more generally \textit{sequipowers}\cite{zimin}. It is then abstractly similar to the celebrated fractal $ABACABA$ sequence\cite{abac}.

The simplest approach to obtain a concrete formula for computing the TZS utilises the binary divisibility by successive powers of two so that one can write\\

$t(x) = \sum_{i=0}^{l_2(x)} 1_\chi(mod(x,2^i)=0)$\\

where $l_2(x) = 1 + \lfloor log_2(x) \rfloor$ is the binary logarithm standing for the maximal power of two of the expansion of $x$. Similar expression will
also arise in higher alphabets in the more general setting posed by Conway's fractran. 
In the case of the binary alphabet, it is also possible to rewrite the same using the Hamming distance between $x$ and $x-1$ as prescribed in the relevant 
OEIS page\cite{oeis}. A graphical representation of the $2^{t(x)}, x\in[0,...,]$ is shown in figure 1.

\begin{figure}[H]
  \centering
  \includegraphics[width=10cm]{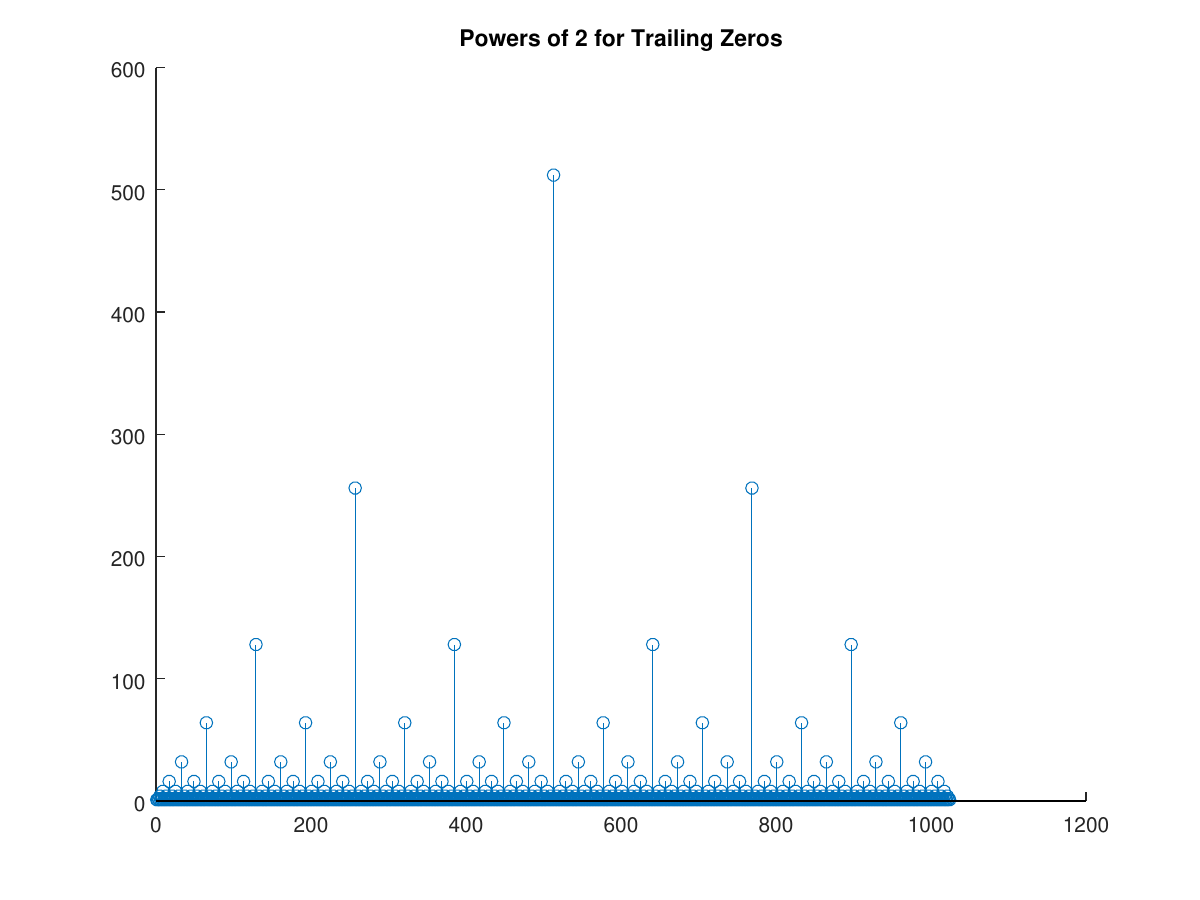}
  \caption{The tree structure of the TZS associated even factors}
\end{figure}

Since, knowledge of a complete shift in the exponents of every 2 factor is possible beforehand, it should also be possible to rephrase the original problem
so that any computation would spend only a single step to each of the two branches given a complete knowledge of the total shift $2^{-t(x)}$ thus effectively
realising the whole bunch of $f_0^{t(x)}$ total application of the first map.

Such a transcription is facilitated by rewriting a composite which takes into account both cases of $f_0^n\circ f_1$ and $f_1\circ f_0^m$ at once by noticing the 
equivalence of $f_1$ with $3x+mod(x,2)$ in which case the orignal CH in (1) is rewritten as

\begin{equation}
x_{n+1} = 2^{-t(x_n)}(3x + mod(x_n,2)) + 1 - mod(x_n,2)
\end{equation}

The expression in (2) automatically sends to either of the two branches depending on the $(mod 2)$ class. It is preferable to rewrite it also in the form\\

$x_{n+1} = A(x_n)x_n + B(x_n)$\\

where $A(x_n) = 3/2^{t(x_n)}$ and $B(x_n) = 1 + mod(x_n, 2)( 2^{-t(x_n)} - 1 )$. 

It is immediately obvious that $B(x_n) = 1$ for all cases. This is simply because of the complementarity of $t(x)$ with $mod(x, 2)$ since all roots of $t(x)$ 
are odd integers. Therefore the final reduction of CH in (1) is equivalent to the expression

\begin{equation}
x_{n+1} = \left(\frac{3}{2^{t(x_n)}}\right)x_n + 1 
\end{equation}

This final form will exactly reproduce the same elements of the original iterations that do not contain successive binary shifts. A possible termination condition for this type of iteration can be given as
$mod( log2(x_n, 2), 1) = 0$. Checking 

The coefficient that apepars in (3) is of special importance and stands for the bridge between the original problem and that of the palindromic binary strings
as explained in the next section. The fixed points of the final map in (3) are found via the standard condition $f(x)-x=0$ rewritten as

\begin{equation}
\left ( 1 - \frac{1}{x}\right )2^{t(x)} = 3
\end{equation}

Given the strucure of the TZS for any expansion of length $L$ in a maximal $L$ interval $[0,...,2^L-1]$ it holds that $0 \le t(x)\le L-1$.
Restricting search in all powers of 2, ($x = 2^l : t(x) = l$) immediately cancels out the exponential term in (3) leaving only the condition
$2^l - 1 = 3$ so that the only possible integer root in (4) is at 
$x = 4$. 

In the next section some appropriate definitions are introduced which will make it possible to establish the connections of the
new single branch map with the issue of palindromic words in fixed length binary expansions.

\section{ Hierarchies of Palindromes}

\subsection{Preliminary definitions}

The particular construct presented requires the introduction of constant length binary expansions for all words inside an interval. When expressed this way, all binary patterns inside an exponential
interval are said to form a so called, \textit{'Hamming Space'} 
the reasoning being that all such expansion are then forming a normed, linear vector space the norm being given be the Hamming distance\cite{ham}. 

All such representations require that any binary expansions are also
characterized by a number of leading zeros. This leads to an additional ambiguity with the definition of certain operators acting on words like reflections or mirror inversions and the construction of palinrdomic words due to the need for additional parametrization for the expansion length required. 

Because of this necessary to setup a different than usual representation which is only possible across a self-similar hierarchy of lexicographically ordered sets that can be represented as special asymmetric matrices of all patterns. To do this, the following terminology will be useful. 

A number $M(L)=2^L-1$ is to be called a \textit{Mersenne} number and an interval $s(L)=[0,...,M(L)]$ is to be called a 
\textit{Mersenne interval}. A self-similar sequence of intervals 

$$
s(1)\subset s(2) \subset\cdots\subset s(L)\subset\cdots
$$

is to be associated 
with a set of $L\times 2^L$ matrices of lexicographically ordered bit patterns as a representation of each $S_L$ in $1-1$ correspondence 
with the binary expansion of the row index $j\in s(L)$ via the polynomial representation.

The particular choice is justified by a variety of reasons including the
fact that the above is also a well formed hierarchy of closures for 
certain binary operators like the \textit{bitwise XOR} which is known
to have the group property.

It is also a known fact that each column of any $S_L$ matrix representation is identical with the paths of a symmetric, homogeneous rooted binary tree thus corresponding to a self-similar
hierarchy of binary tree structures. 

Equivalently, the same can be
phrased as an arithmetic equivalent of a hierarchy of \textit{Hamming Cubes} or, subspaces of an $L$-dimensional hypercube, due to the fact
that every element of a Hamming space associated with a fixed length
binary expansion can be put into a one-to-one association with the 
edges of such a hypercube\cite{cubes}

The particular form of the hierarchy of lex-ordered 
matrices is also known in another context as a set of \textit{Orthogonal Designs}\cite{lex} when written in the equivalent $\{\pm 1\}$ alphabet 
instead of $\{0, 1\}$.

To further facilitate an exchange between the language of sequences and binary patterns of constant length, it is appropriate to denote 
$|w|\in \mathbb{N}$ for the arithmetic value of each binary word via the use of an "encoding" map

$$
p: |w|=p(w) = \sum_{i=0}^{L-1}a_i2^i
$$

and its abstract "decoding" inverse 

$$
p^{-1}:[a_0,...,a_L] \overset{ p^{-1}}{\leftarrow}|w|
$$

This helps establishing a direct association of the hierarchy of intervals with the hierarchy of matrices as\\

$$
\begin{tikzcd}[cells={nodes={minimum height=2em}}]
s(0) \arrow[r, "\delta"]  & s(1) \arrow[r,"\delta"] \arrow[d, "p^{-1}"] & s(2) \arrow[d, "p^{-1}"] \cdots\\  
\left[0\right] \arrow[r] & S_1 \arrow[r, "\Delta"] & S_2 \cdots
\end{tikzcd}
$$

The additive action of the $\delta$ map is to simply increase the cardinality of any previous ordered list representing any Mersenne interval by applying
to each member the same rule and perform a list concatenation like\\ 

$$
\mathfrak{\delta}s(L) \rightarrow s(L+1) = [s(L), s(L) + 2^L], L=0,1,\cdots
$$ 

Correspondingly, each new application of the decoder $P^{-1}$ results in an equivalent self-similar action denoted by $\Delta$ comprising 
a concatenation of a copy the previous $S(L)$ matrix followed by the addition of a new top row of precisely $L$ zeros to be followed by 
a left to right flip of its not-complement like

$$
S_1 =
 \begin{bmatrix}
 0 & 1 \\
\end{bmatrix}
\rightarrow 
\begin{bmatrix}
 0 & 0 \\ 
 0 & 1 \\
\end{bmatrix}
\rightarrow 
S_2 =
\begin{bmatrix}
 0 & 0 & 1 & 1 \\
 0 & 1 & 0 & 1 \\
\end{bmatrix}
$$
 
 Thus, the whole hierarchy comprises a sequence of internal 
 2-complements and reflections.
 
 Effectively, all rows of any member matrix of the hierarchy is a periodic sequence spanning an exponential sequence of periods that can be directly  computed via either a Boolean or an equivalent  arithmetic formula as
 
 \begin{equation}
 S_{i,j\in[0,\cdots,2^L-1]}(L) = 2^{-i}(j\otimes 2^i) = mod\left( \lfloor \frac{j}{2^i} \rfloor, 2 \right), i=0,1,\cdots,L-1
 \end{equation}
 
where $\otimes$ stands for the bitwise AND operation.

 Moreover, one may consider the action of arbitrary automata $\mathcal{A}$ accepting some or all of the rows of each member $S_L$. There are two main classes
 of possible automata described as either a)\textit{Indicators} or maps $\mathcal{A}:S_L\rightarrow \{0,1\}$ asserting existence of 
 some property of a binary  pattern or b)\textit{Transducers}: $\mathcal{A}:S_L\rightarrow S_N$. 
 
 Restricting attention to such automata or their equivalent Turing machine expressions
 that halt, one may introduce an upper bound as $max(L,N)$ such that every such halting automaton represents an endomorphism in $S_{Lmax}$. It is then possible
 to project the action of all such automata into a new sequence or production via 
 
 \begin{equation}
 s_{\mathcal{A}} \leftarrow \left( p\circ\mathcal{A}\circ p^{-1}\right)s(L)
 \end{equation}
 
 This implies the reduction of any such computation into a chain of maps producing a hierarchy of sequences of exponentially increasing length \\
 
 \begin{tikzcd}[cells={nodes={minimum height=2em}}]
s(0) \arrow[r, "\delta"]  & s(1) \arrow[r,"\delta"] \arrow[d, "p^{-1}"] & s(2) \arrow[r,"\delta"] \arrow[d, "p^{-1}"] & \cdots\\  
\left[0\right] \arrow[r] & \mathcal{A}[S_1] \arrow[r, "\Delta"] \arrow[d, "p"] & \mathcal{A}[S_2] \arrow[r, "\Delta"] \arrow[d, "p"] & \cdots\\
s_{\mathcal{A}}(0) \arrow[r,"g"] & s_{\mathcal{A}}(s(1)) \arrow[r, "g"] &  s_{\mathcal{A}}(s(2)) \arrow[r,"g"] & \cdots\\
\end{tikzcd}
\\
 
 Computability of the action of the new $g$ map in terms of simple 
 arithmetic formulas is a separate and difficult issue in gerneral for arbitrary automata. Whenever possible, one may exchange the action of
 any such automaton with the resulting sequence via a recusrsive list
 concatenation or even the generating function of the resulting sequence
 if summable at all.
 
 A trivial example can be given in case of an automaton of which the 
 action equals some permutation $\pi$ of symbols across $w\in S_L$ 
 in which case one can apply a 'transfer'
 principle based on the constant expansion length as
 
 $$
 \sum_{i=0}^{L-1}a_{\pi(i)} 2^i \cong \sum_{i=0}^{L-1}a_i 2^{\pi(i)} 
 $$
 
 For instance, inversion of all positions resulting in a
 mirror inversion of all words will also result in a map of the original
 sequence of natural numbers in an iterative sequence given by
 
 $$
 r(i+1) \leftarrow [r(i), r(i) + 2^{L-i-1}], r(0) = 0
 $$
 
 Another simpler alternative of scaling maps for reflections is presented in the next section. 
 
 Notably, the sequence of applications of any such scaling map follows 
 a well known pattern of another fundamental fractal sequence the so
 called, \textit{Sum-of-Digits} sequence\cite{sd} as 
 
 $$
 I\, g\, g\, g^{(2)}\, g\, g^{(2)}\, g^{(2)}\, g^{(3)},\cdots,g^{sd(n)}
 $$.
 
 following the pattern
 
 $$
 0\ 1\ 1\ 2\ 1\ 2\ 2\ 3\ \cdots
 $$
 
 The $sd(n)$ sequence also admits the simplest arithmetic scaling map $g(x)=x+1$
 following the linear staircase of the natural lengths or maximal powers
 of two present in any lexicographically ordered set of binary patterns.
 
 Sequences for which the leading zeros may play a role in their definition will not admit as simple a recursion as above and it will
 often exhibit a similar recursive structure with a branched map acting
 differently on the first and second part of a list concatenation scheme
 in the abstract form
 
 $$
 g(s(i+1))\leftarrow [(g_0)(s(i),i), (g_1)(s(i), i)] 
 $$
 
 where the iteration index $i$ may have to be explicitly included in general. Equivalently, the resulting compositions can always be
 extracted by a symbolic binomial expansion $(g_0 + g_1)^L$. 
 The particular case of fixed length reflections and palindromes is  analyzed in the next section.

\subsection{Palindromes and the fixed length reflection group}

Let $w$ a binary word and $\mathcal{R}(w, L):S_L \rightarrow S_L$ also denoted as $R_L$ heretofore, an order reversing map also called the \textit{reflector} heretofre, with $S_L$ the set of all $2^L$ binary strings of same length $L$. 

Then if $w = [a_0,\cdots,a_L]$, its mirror inversion or reflection is denoted as $\mathcal{R}_L(w) = [a_L,\cdots,a_0]$. 
$\mathcal{R}_L$ is then one of the two natural 
involutions for any words in every $S_L$ the other being the 2-complement.

The hierarchical construction of the previous section
imposes a discrimination between even and odd order palindromes depending on $L$ being even or odd as well. Thus for all odd order matrices $S(2k+1)$
a palindrome may leave the "central" symbol at $k+1$ unaltered and only invert the order of the first $k$ symbols so that 
\\

$a_{2k+1} = a_0,...,a_{2k+1-i}=a_i,...,a_{k+2} = a_k, i=0,1,...,k$\\

This class will not be treated further heretofore.\

While the 2-complement is fixed point free in any $S_L$, the reflection
operation will admit a set of fixed points known as palindromic words or
simply, palindromes. 

By construction the number of fixed points of every even order $S_{2k}$ must
necessarilly contain the whole reflected set in $S_k$ since anyone of them
can get reflected thus forming a member of the set of fixed points 
of the reflecor operator over any $S_{2k}$. Hence, the cardinality of
the sets of fixed points must also form a sequence of $2^k$ fixed points
of the $\mathcal{R}_{2k}$ action over any $S_{2k}$ across the hierarchy. Obviously the 'edges' of each set comprising the all
zeros and all ones pattersn are always fixed points. 

From now on the notation $\mathcal{R}_L(|w|)$ will be interpreted as the 
expanded form of (6), that is $\left( p\circ\mathcal{R}\circ p^{-1}\right)(|w|)$ which acts on the valuation of the word $w$ by expanding, processing and
again contracting to a new integer. Thus the total action over the sequence
of natural numbers inside any Mersenne interval wiil always result in a
new sequence in that same interval parametrized by the additional length
parameter $L$.

By definition, a reflector has the group property since it sends one to
one, any integer inside the same Mersenne interval thus being equivalent to a permutation. This is one of the main reasons for using the hierarchy
over fixedl length expansions. Otherwise, eny power of two would be
mapped to a one after bit order reflection thus failing to be a bijection.

The reflector does not act homomorphically over the standard arithmetic
addition and multiplication but it does obey an antihomomorphism with 
respect to the arithmetic equivalent of concatenation\\

$\mathcal{R}_L(|x|+|y|2^L) = \mathcal{R}_L(|y|) + \mathcal{R}_L(|x|)2^L$\\
 
The difference of fixed length reflections with leading zeros becomes more evident by noticing the appearence of a bit shift 
in any newly formed subsequence of reflected integers due to
inversion of position of the leading zeros blocks. 

Thus, a fixed length reflection can also be subsumed via the use of two scaling maps corresponding to the $g$ map of the chain diagram of previous section leading to two different iterative list concatenation methods as\\

$r(s(L+1)) \leftarrow [g_0r(s(L)), g_1r(s(L))], r(s(0)) = 0$\\

where now $g_i(x) = 2x + mod(i,2)$ applied pointwise across all
previous list elements.

Before coming in the subject of fixed length palindromes and 
their properties it is important to add another toolbox in the description of fixed length binary expansion property of every
fixed length reflected binary expansion which binds them with 
the TZS as well as the leading zeros sequence or LZS in a
particularly useful way.

For all such binary words, an alternative compressed
representation exists given in terms of a \textit{run-length}
encoding in the form of an alternating polynomial given as a
bijective map\\

$rl(w):S_L\leftrightarrow \mathbb{Z}^L:[a_0,\cdots ,a_{L-1}]\leftrightarrow [\pm c_0(|w|)\cdots ,\pm c_m(|w|)]$\\

under the convention of a minus sign for a block of zeros and vice versa.
Each coefficient $c_i$ counts the length of a block of same
symbols marking with a $\mp$ sign whether it is a zero or one
respectively.

An additional constraint over all alternating coefficients results
from the fixed length expansions in the form

$$
\sum_{i=0}^{m(n)}|c_i| = L, n \in s(L) 
$$

Due to the bijective nature of this mapping all distinct \textit{integer partitions}\cite{ipart} of $L$ exist inside the
set of RL representations of the members of any $S_L$ since any such can always be turned to a binary pattern of same length. 

The first and last of these coefficients, when expressed as
sequences over all integers in the associated Mersenne interval
are of special importance. In particular, the TZS is equivalent 
to all negative values of the first coefficient at even indices
$c_1(2k)$, and zero for all even indices while the leading 
zeros sequence (LZS) is identified with $c_m(n), n\in s(L)$.

The latter is naturally associated with the binary logarithm
$l_2(n)$ via $c_m(n) = L - l_2(x)$ in accord with the fixed 
maximal length representation used here.

The three fundamental sequences characterizing each binary pattern
given by the triplet ${l_2, t_n, ds_2}$ share the same range in
$[0,\cdots, L]$. Moreover, the TZS and the max. bit sequence
$l_2$ share the same multiplicities of values. 

The following proposition can also be proven

\begin{prop}
\begin{enumerate}
 Let $t_n$ be a sequence of all $t(n),n\in s(L)$ and let $r_n$ the 
 sequence of reflected binary expansions as integers over the same interval. 
 Then for all  $L$, the second is a decreasing order sorting
 permutation of the  first or $(t_z\circ r)(n)\cong sort_{>}(t_z(n))$.
\end{enumerate}
\end{prop}

This is a trivial result of the simultaneous fixed length
reflection of all binary words in any interval which affects 
an exchange of the first and last coefficients in the RL
representation and hence of the TZS with the LZS.

The corresponding permutation then is characterized by the first
blocks of zeros being already sorted in size due ot the fact that
any leading  zeros will have a difference from the maximal power of two that scales as $L - i$ with every new exponential
subinterval $s(i)$ which adds a single bit on top of all the previous expansions.

As a result, the permuted integers will contain a decreasing
number of all even integers with $2^{L-i-1}$ factors when
lexicographically ordered on the first half of the interval with all odd integers also mapped to the second half carrying over
all zero values of the TZS. 

Subsequent omposition with the standard form of the TZS is then
equivalent to a sorted counting of the number of unique digits 
in a TZS sequence of the form\\

$ 1\ 0\ 2\ 0\ 1\ 0\ 3\ 0\ 1\ 0\ 2\ 0\ 1\ 0\cdots $\\

The result of such a composition of sequences is shown in figure 2
on all integers in $[0,\cdots 2^{10}]$ in a semilog graph to
make the staircase structure more pronounced.

\begin{figure}[H]
  \centering
  \includegraphics[width=10cm]{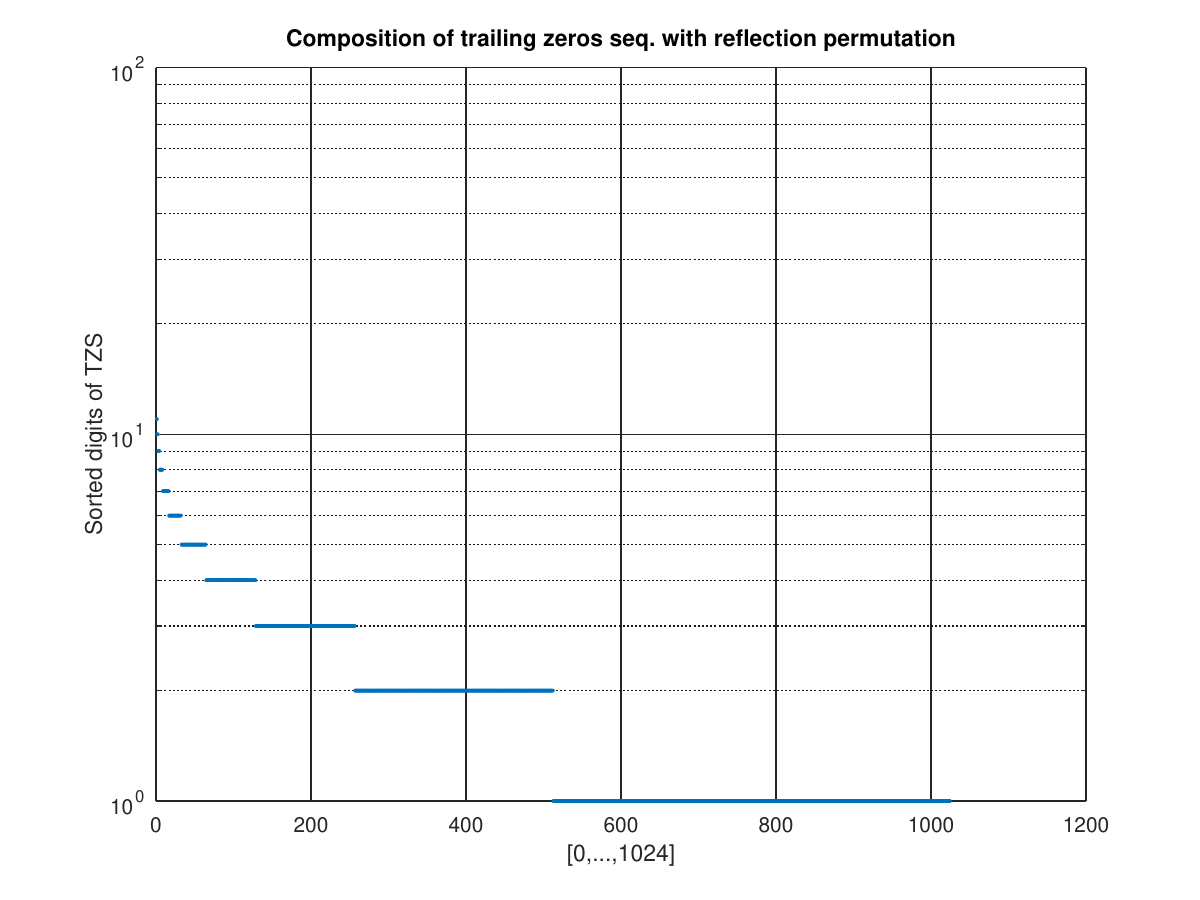}
  \caption{The composite sequence of the TZS over all reflected integers}
\end{figure}

One can then use directly this result to extract an integer 
histogram of all unique digits in any subsequence of TZS over any
Mersenne interval and which follows a scaling law of the 
form $2^{L-l_2(i)-1}$ in accord with the tree structure also shown in figure 1.

Furthermore, it can be inductively verified that for any maximal
interval $s(L)$ the following relation is always satisfied

\begin{equation}
l_2(x) + (t\circ r)(x) = (l_2\circ r)(x) + t(x) = L_{max}
\end{equation}

where $L_{max} = Max(l_2(x))$. 

Given the action of the reflector it is also possible to introduce an arithmetic decomposition of every palindrome's integer value of even order as

\begin{equation}
\mathcal{P}_{2k}(|w|, L) = \mathcal{R}_k(|w|) + |w|2^k
\end{equation}

The expression in (8) utilizes the involutory nature of the reflector
so as to get a sorted sequence of all possible even palindromes in
any $S_L$. This is simply the result of taking two copies of any 
$S_L$ matrix and perform a horizontal concatenation of the first copy
with an up-down flip of the second copy, something easy to realize in
an array language like Matlab or Octave.

Getting back the corresponding integer values will always result 
in a sorted sequence since the single application 
of $\mathcal{R}_L$ is equivalent to a permutation but the addition 
of all original bits above $2^L$ guarantees this being a sorted
sequence.

Let then, $\mathfrak{\Delta}\mathcal{P}_L$ denote the discrete
differences of the sorted sequence in (8) being again decomposable as

\begin{equation}
\mathfrak{\Delta}\mathcal{P}_L(|w|, L) = \mathcal{P}_L(|w|, L) - \mathcal{P}_L(|w|-1, L) = \mathfrak{\Delta}\mathcal{R}_L(|w|) + 2^L
\end{equation}

where now

\begin{equation}
 \mathfrak{\Delta}\mathcal{R}_L(|w|) =  \mathcal{R}_L(|w|) - \mathcal{R}_L(|w|-1)
\end{equation}

It is then possible to prove the below proposition

\begin{prop}
\begin{enumerate}
for all $|w|\in \mathcal{M}(k)$ and their expansions $w \in S_{k}$ for any $k$ it holds that \\

$\mathfrak{\Delta}\mathcal{P}_{2k}(|w|) = 3(2^{k-t(|w|)-1})$\\ 

where $t(|w|)$ the trailing zeros sequence.
\end{enumerate}
\end{prop}

The proof can be given with the aid of elementary curry-less binary addition
performed on an unbounded or bi-infinite tape of a TM restricted in the case
of a successor function $\mathfrak{s}(|w|) = |w|+1$.

Indeed, adding a single bit at the lower power of any expansion only requires
two rules. Assuming any integer $n$ coded in binary with powers of two from
left to right, these rules are

\begin{itemize}
 \item if $n$ is even, the head writes '1' in the present position and stops.
 \item if $n$ is odd, the head moves to the left replacing all 1s with 0s until it reaches the first 0 position where it writes a '1' and stops.  
\end{itemize}

Next, consider the case of a bi-infinite tape with the first digit situated
at a central cell with the whole pattern reflected as if by a mirror in the 
middle. One can alwasy assume a machine with two heads working in opposite
directions but following the same pair of rules.

There can only be transitions from even to odd or from odd to even
numbers. In both cases, any alterations in the first block of digits will
not affect the next blocks so that they cannot contribute to the discrete
diffferences of the resulting sequence of palindromes.

In the first, even-to-odd case, assume an arbitrarily large all 0s block
reflected across the middle point. Then the transition will be of the form

\begin{itemize}
 \item $a_0\ \cdots 1\ 0\ \cdots 0 | 0 \cdots 0\ 1 \cdots\ a_{2k} $
 \item $a_0\ \cdots 1\ 0\ \cdots 1 | 1 \cdots 0\ 1 \cdots\ a_{2k} $
\end{itemize}

Since the whole palindrome is now a new pattern with an expansion of
double length, the newly added 1s in the middle will correspond to a
pair of new powers $\{a, 2a\}$ with $a=2^{k-1}$. Then inevitably, 
the difference between the new integer advanced by one will have to be 
$3a = 3\ 2^{k-1}$.

In the second, odd-to-even case, assume again an arbitrarily large all 1s
block in which case any transition will be of the form

\begin{itemize}
 \item $a_0\ \cdots 0\ 1\ \cdots 1 | 1 \cdots 1\ 0 \cdots\ a_{2k} $
 \item $a_0\ \cdots 1\ 0\ \cdots 0 | 0 \cdots 0\ 1 \cdots\ a_{2k} $
\end{itemize}

By a similar argument as before, the new 0s block must contain a
number of $2t(|w|)$ 0s. With both blocks shifted by the same amount of
$a = 2^{k-t(|w|)-1}$ any difference becomes\\

$(2^{2t+1} + 1)a - 2(2^{2t} - 1)a = 3\ 2^{k-t(|w|)-1}$

Consequently, we also obtain the following computationally useful results.

\begin{corollary}

 A. The sequence of valuations of palindromic words over any $s(k)$ interval is given by the sequence of partial summands\\
 
 $\mathcal{P}_{2k}(i) = 2^{k-1}3\sum_{i=1}^{2^k}2^{-t(i)},$\\
 
 B. The sequence of reflectors over any $s(k)$ interval is given as\\
 
 $\mathcal{R}_k(i) = 2^{k}\left( 3\sum_{i=1}^{2^k}2^{-t(i)-1} - 1\right)$
 
\end{corollary}

In the light of the relation (7) the previous can be further
simplified as

\begin{eqnarray}
 \mathcal{P}_{2k}(i) & = &\frac{3}{2}\sum_{i=1}^{2^k}2^{-(l_2(r(i))}\\
 \mathcal{R}_k(i) & = &\frac{3}{2}\sum_{i=1}^{2^k}2^{-(l_2(r(i))} - 2^k 
\end{eqnarray}

One immediately notices here the presence of the magic factor
$3/2^{t(x)}$ in the total expression of the corollary $1$ for 
$\mathfrak{\Delta}\mathcal{P}$.
This is then used to redefine the dynamics of the modified CH in 
(3) in a particular way revealing an interaction between two 
'mirror' worlds.

\section{A hidden mirror in the CH dynamics}

It is obvious from direct comparison of the expression of the map in (3)
and the result in proposition $2$ that one should be able to make a direct substitution as

\begin{equation}
x_{n+1} = \left(\frac{\mathfrak{\Delta}\mathcal{P}_{l_{max}}(x_n)}{2^{l_{max}-1}}\right)x_n + 1 
\end{equation}

where now the coefficient numerator in (10) is to be interpreted as

\begin{equation}
\mathfrak{\Delta}\mathcal{P}_{l_{max}}(x_n) = \mathcal{P}_{l_{max}}(x_n) - \mathcal{P}_{l_{max}}(x_n-1) = \mathfrak{\Delta}\mathcal{R}_{l_{max}}(x_n) + 2^{l_{max}}
\end{equation}

In order for the substitution to make sense it has to be assumed that
each $x_n$ is taken inside an interval $s(l_2(x_n))$ which varies.

For this reason it is necessary to use a varying maximal length for
the definition of the reflector as $l_{max} = Max(l_2(x))$ (practically $2^{l_{max}}$ is equivalent to the use of standard libraries like \textit{nextpow2}). This becomes necessary due to
the fact that there is a hierarchy of different $\mathcal{R}$ sequences or reflective permutation groups across different intervals $s(L)$.
On the other hand, the original conjecture is equivalent to the
existence of an upper bound for all such intervals.

The dynamics in (13) appears now as the result of an interaction of
the original variable with the 'slope' formed by a discrete
derivative over reflections. To further understand this version of
the original dynamics it is necessary to find a reduction of the
discrete difference in some more fundamental sequences like those
introduced in the previous section. 

\begin{figure}[H]
  \centering
  \includegraphics[width=10cm]{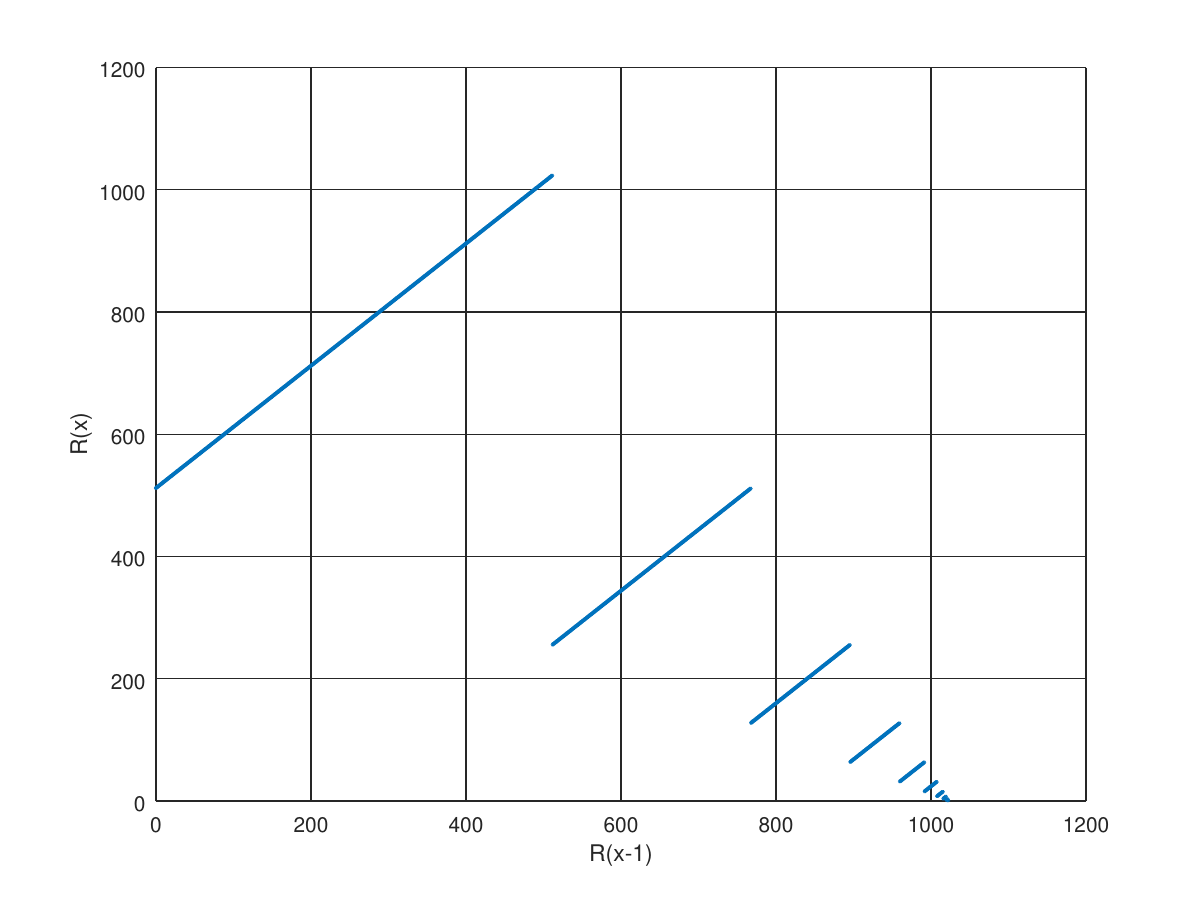}
  \caption{Global map for the reflection sequences}  
\end{figure}

The simplest way is to utilize a global map for the pair 
$\{\mathcal{R}(x), \mathcal{R}(x-1)\}$. This is
shown in figure 3, where a scaling law appears to govern the 
piece-wise linear dependence over successive intervals.

An arithmetic formula can be found inductively across the hierarchy
of intervals utilizing the particular scaling which follows a 
similar pattern with that of the sorted TZS in figure 2. 

It is then possible to prove inductively for the map of figure 3, 
the arithmetic interpolant

$$
 \mathcal{R}_L(x) = \mathcal{R}_L(x-1) - 2^{l_{max}} + 3\times2^{\mathfrak{l}(x)-1}
$$

where now

$$
\mathfrak{l}(x) = l_2(2^{l_{max}} - x - 1)
$$

It should be noticed that the argument in $\mathfrak{g}$ 
performs a kind of parity reflection over any interval due to the 
arithmetic equivalent of the 2-complement defined as $NOT(x,L) =
2^L - 1 -x$. 

Then the original reflector difference reduces to elementary
sequences as

\begin{equation}
 \mathfrak{\Delta}\mathcal{R}_{l_{max}} = 2^{l_{max}}\left(  3\times2^{\mathfrak{l}(x)-l_{max} - 1} - 1\right)
\end{equation}

Substitution in (13) using (14) then results in 

\begin{equation}
x_{n+1} = 2\left(\frac{\mathfrak{\Delta}\mathcal{R}_{l_{max}}(x_n)}{2^{l_{max}}} + 1\right)x_n + 1 = 2^{\mathfrak{l}(x_n)-l_{max}}(3x_n) + 1
\end{equation}

What is actually gained in (16) is the expression of the same 
dynamics as in (3) but this time avoiding the difficulty of the 
fractal structure of the TZS sequence. 

On the other hand, in the light of relation (7) it is also possible
to write (3) as

\begin{equation}
 x_{n+1} = 2^{l_{max}-l_2(\mathcal{R}(x_n))}(3x_n) + 1
\end{equation}

The expression in (17) again emphasizes the role of reflections in
the overall dynamics evident in the antagonism between the
effective lengths of two mirror images in the exponent of (17). 

Actually, the two expressions have now come full circle since
the parity reflection in $\mathfrak{l}$ is just another identity
in disguise or $\mathfrak{l}(x)-(l_2\circ r)(x) = 2l_{max}$ on 
the exchange of parity reflections with the index binary reflections. 

This alone is not sufficient to explain the mystery of the 
conjectured global convergence yet another alternative is offered
in the last section which may be fruitful for further investigation
in juxtaposition with the type of 'mirror' image interactions
presented.

\subsection{Convergence as block decimation}

From the structure of (3) it is evident that any final convergence
to the fixed point of the dynamics will take place as soon as the 
trajectory will reach a pure power of two.

In the light of the equivalent RL representation introduced in
section 3.2, this can be phrased as a reduction of the number
of blocks since any such is always of the form $\{\mathbf{0}^{c_0},\mathbf{1},\mathbf{0}^{c_2}\}$ for any $l_{max}=c_0+c_2+1$.

The originally conjectured global convergence must then be equivalent
to a higher probability of a falling number of blocks leading to what
could be termed a \textit{'block decimation'} effect although it is
not stepwise homogeneous. A possible strategy for proving the
original conjecture could then start with a proper definition of 
such a probability.

An effective measure of such a type of binary complexity index
can be given in terms of the number of RL coefficients which may be
termed here as the \textit{RL Dimension} or RLD for brevity. 

The particular type of RLD sequences per $s(L)$ interval also have a
fractal character and can be found via induction over the hierarchy 
to satisfy a scaling law given by the standard list concatenation
scheme
$$
 rld(i+1) \leftarrow [rld(i), \mathcal{R}(rld(i)) + \sigma_i], rld(0) = 1
 $$
 
 where now $\sigma_i = 1,i=0,\cdots L-1$ and $\sigma_L=0$ while
 the reflection operator inverts the list index order at every step.
 
 A natural property of any such sequence is its invariance under
 the reflection group over the indices themselves or simply 
 $rld(\mathcal{R}(n)) = rld(n),n\in s(L)$ since reflection over 
 fixed length  binary expansions of each index cannot alter the
 number of the corresponding RL coefficients.
 
 An example of such a fractal sequence can be seen in figure 4(a)
 over $s(10)$ while in figure 4(b) its distribution is compared
 against the standard binomial distribution of the 'Digit-Sum'
 sequence.
 
\begin{figure}[H]
  \begin{subfigure}[b]{0.6\textwidth}
    \includegraphics[width=\textwidth]{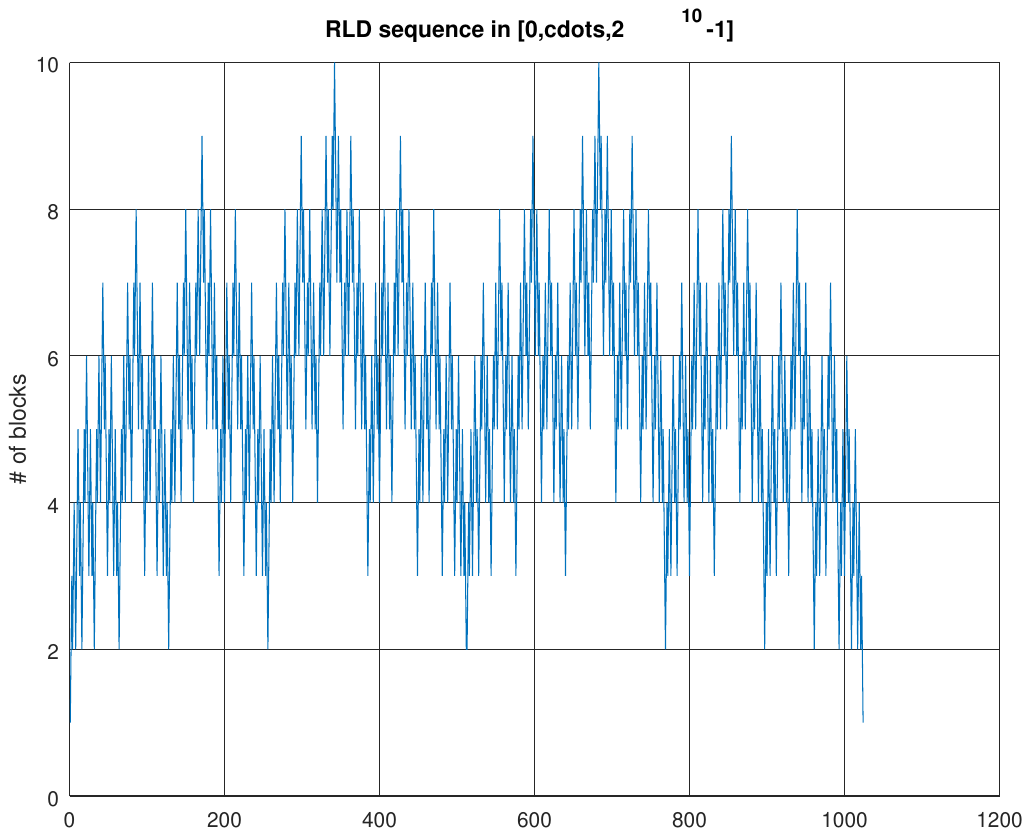}
  \end{subfigure}
  \hfill
  \begin{subfigure}[b]{0.6\textwidth}
    \includegraphics[width=\textwidth]{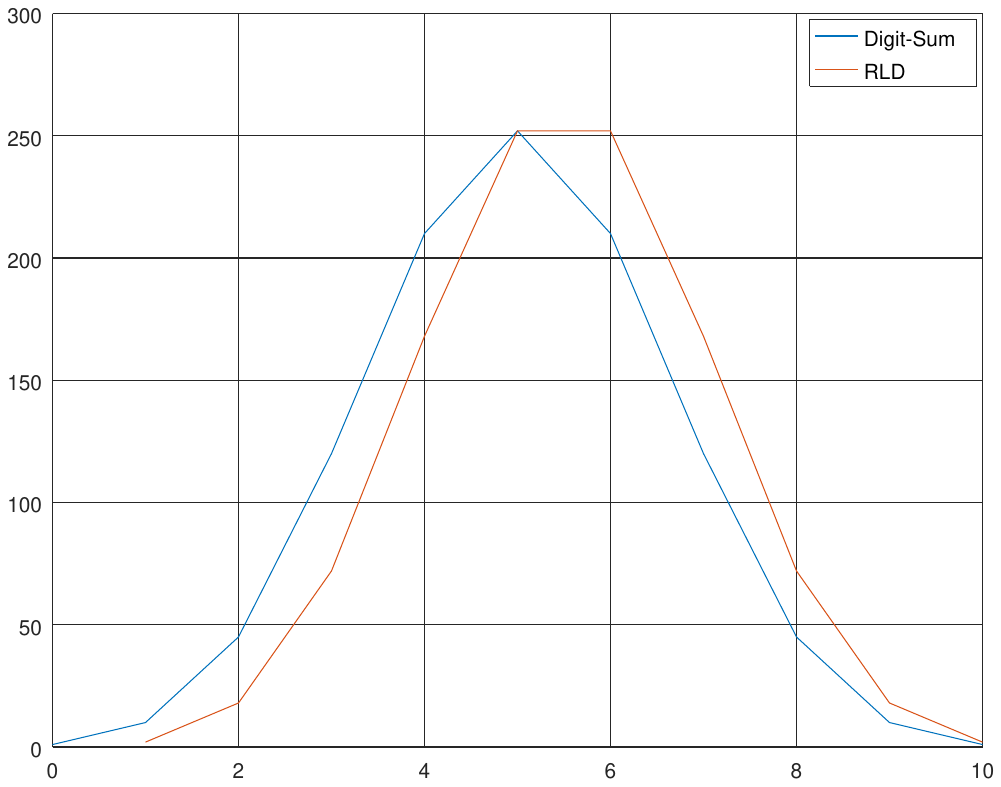}
  \end{subfigure}
  \caption{ (a) Example of an RLD sequence on $s(10)$ and (b) Distributions for 'Digit-Sum' and RLD on $s(10)$ }
\end{figure}

 The simplest way to utilize these sequences is to define a 
probability of any odd integer being mapped to another even one
via the '$3x+1$' map of smaller or larger RLD. 
Given a long list of RLD values, this
is straightforwardly written as

\begin{equation}
 p_{<,>} = \frac{1}{\mu_L(rld)}\sum_{x=2k+1,x\in s(L)} sign\left[rld(3x+1) - rld(x)\right] 
\end{equation}

where $x$ runs in all odd values in $s(L)$ and $\mu_L(rld)$ is an
appropriate normalization measure over all $2^L$ values of RLD.

A preferable index for asserting any increase or decrease in the 
complexity of the resulting patterns can then be given by the 
probability mass ratio $p_</p_>$ which avoids normalization.

Additionally, the previous section finding suggests a correlation
of such an index with the 'interaction' between mirror images
of the expansions of the $x_n$ variable via the quantity

\begin{equation}
\delta_R(x) = l_2(3x+1)-l_2(\mathcal{R}(3x+1)) 
\end{equation}

The associated probabilities 

\begin{equation}
q_{<,>} = \frac{1}{\mu_L}\sum_{x=2k+1,x\in s(L)}sign(\delta_R)(x) 
\end{equation}

allow defining another mass ratio as $q_</q_>$. 

The new ratio can be used to check whether there is 
an increase or decrease in the resulting effective expansion length
between these images. A decreasing ratio could be associated with 
a possible increase in the presence of large binary shifts thus
diminishing the number of blocks to be eliminated.

Indeed, numerical evidence is in favor of this assumption. Results
for a set of itnervals from $s(4)$ up to $s(24)$ are shown in 
figure 5 for both mass ratios. 

Interestingly, the legnths ratio appears saturated soon after the 
tenth power of two while the block ratio increases almost
constantly in a log scale.

\begin{figure}[H]
  \centering
  \includegraphics[width=10cm]{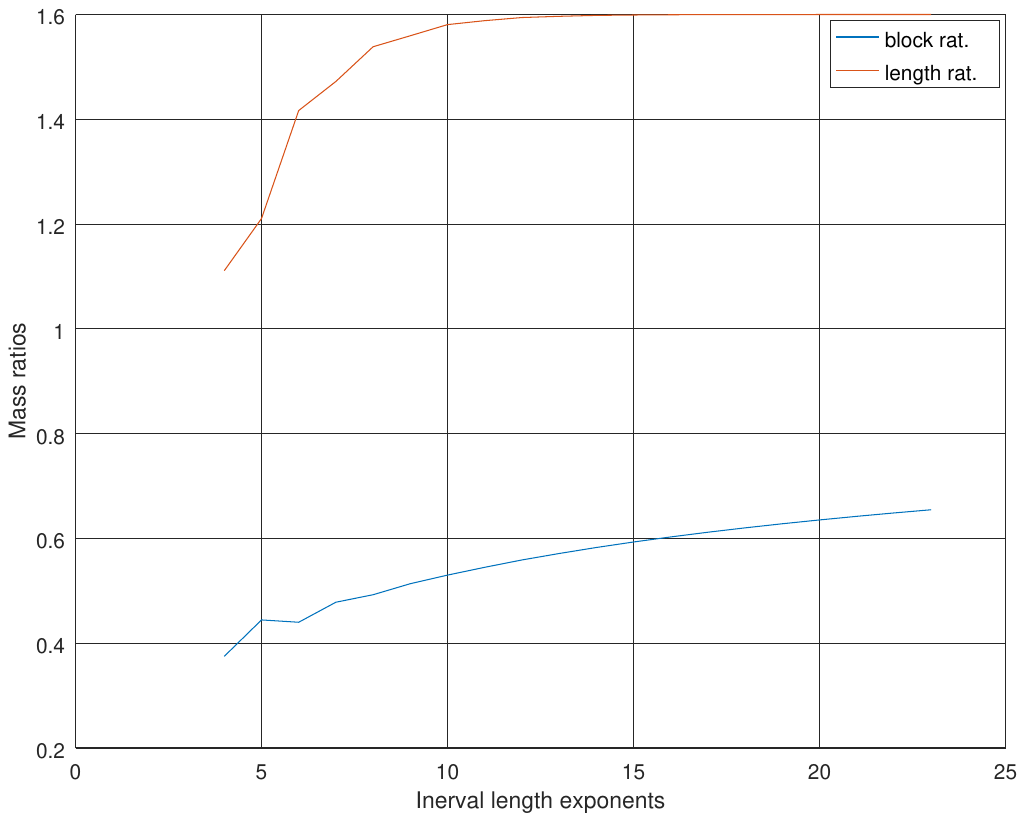}
  \caption{Combined probability mass ratios}  
\end{figure}

\section{ Discussion and Conclusions }

A methodology for combinatorics of automata was introduced which
may offer certain advantages regarding the extraction of scaling
maps and recursive relations over lists as representatives of 
properties of fixed length binary patterns.

The particular application in the case of the Collatz-Hailstone
dynamics was based on a coincidence after reformulating the 
original problem in a way that naturally incorporated the 
origina lbranching condition via tha use of a number theoretic
function known as the trailing zeros sequence (TZS) otherwise
known as the 2-adic valuation of the integers.

When comparing the new form with the sequence of discrete 
differences of palindromes defined via tha action of the 
rteflection group on a hierarchy of exponential intervals 
or closures over the integers, a deeper relation was recognised 
and further analysed.

It was revealed that internal reflections of the binary forms
hidden behind the production of Hailstone sequences play a role
not yet well understood. From a physicist's perspective
it is tempting to think of this dynamics as a bistable potential
with a middle barrier separating two mirror worlds perhaps
amenable to noisy perturbations. Furthering this treatment is
reserved for a future report.

Additionally, there are still unexplored issues regarding the 
multivalued character of the TZS. As a matter of fact, the 
particular symbolic substitution used in (13) of section $4$
could be generalized by simply allowing the index of the 
palindromic word or its internally contained reflection to
be associated with any integer pre-image of the TZS giving
the same value say as $\mathcal{R}(x_n)\rightarrow \mathcal{R}(x_n)\pm \sigma_n$
where $r$ a random variable restricted each time to the same 
level of the TZS.

This brings about an interesting association of the natural 
tree structure of the TZS with a well known mechanical analog
of the so called, '\textit{Quincunx}' or \textit{'Galton Machine'}
\cite{galton} One can think of the additional random variable
$\sigma_n$ as the result of a falling ball across the quincunx
board made out of integerl spacings with its associated height variable thus recreating the exact same values of the TZS.

This leads to the amazing observation that despite $\sigma_n$
being binomially distributed the Hailstone sequences would 
remain absolutely insensitive and hence the dynamics of this 
type would be an invariant of such a perturbation!

It is an ambitious project to carry over similar generalizations
that overcomes the scope of the present short report which was
based on a rather trivial original observation yet it was laid
here in the hope that it may be of aid in future attempts towards
a formal proof of the original conjecture by Collatz.

\end{document}